\documentclass[12pt]{article}

\usepackage{amsfonts}
\usepackage{amsthm}
\usepackage[mathscr]{eucal}
\usepackage{amsmath}
\usepackage{amssymb}
\usepackage{amscd}

\theoremstyle{plain}
\newtheorem{lemma}{Lemma}[section]
\newtheorem{theorem}[lemma]{Theorem}
\newtheorem{proposition}[lemma]{Proposition}
\newtheorem{corollary}[lemma]{Corollary}

\theoremstyle{definition}
\newtheorem{example}[lemma]{Example}
\newtheorem{definition}[lemma]{Definition}
\newtheorem{remark}[lemma]{Remark}

\numberwithin{equation}{section}
\def\begeq{\stepcounter{lemma}\begin{equation}}

\date{}

\begin{document}

\title{The Cone of Semisimple Monoids\\ with the same Factorial Hull}

\author{Lex E. Renner}
\date{February 2006}
\maketitle


\begin{abstract}

\noindent The {\em factorial hull} of the projective variety $X$ (or its cone) is a 
graded algebra $R(X)$ that can be used in some situations to consider simultaneously all divisor classes on $X$.
In this paper we consider initially the situation where $X$ is a {\em semisimple variety} associated with the
semisimple monoid $M$. The factorial hull of such $X$ is determined by a certain arrangement $\mathbb{H}$ of
hyperplanes in the space of rational characters $X(T)\otimes\mathbb{Q}$ of a maximal torus $T$ of $G_0$. If
$G_0$ is simply connected $R(X)$ is the coordinate ring of Vinberg's enveloping monoid $Env(G_0)$. Associated
with $X$ is a certain cone $\mathscr{H}\subseteq Cl(X)$ in the class group
of $X$.  Each $\delta\in\mathscr{H}$ corrresponds to a semisimple monoid
$M_\delta$ with $R(M_\delta)=R(X)$. $M$ and $N$ have the same factorial hull if $X_M$ and $X_N$ differ by $G\times
G$-orbits of codimension two or more. We calculate $\mathscr{H}$ explicitly in the case where $X$ is the
{\em plongement magnifique} for the simple goup $G_0$. This is exactly the
case where $M$ is a {\em canonical monoid}.

\end{abstract}


\section{The $\mathbb{Q}$-Factorial Hull}     \label{one.sec}


The $\mathbb{Q}$-factorial hull of a projective variety is easy to describe, assuming it exists.
This construction is related to a well-known question of Hilbert \cite{Mum,Nag,Rees,Zar}.
Assume that $X$ is an irreducible, normal, projective variety over the algebraically
closed field $K$ of characteristic zero. Let $Cl(X)$ be the divisor class group of $X$, and assume that
$F\subseteq Cl(X)$ is a free abelian subgroup of finite rank.
Choose representatives $M_\alpha$, $\alpha\in F$, consisting of rank-one,
locally reflexive sheaves on $X$. By results of \cite{Cox,Hausen,E-K-W}, there is a
natural graded $K$-algebra structure on 
\[
R_F(X)=\bigoplus_{\alpha\in\Gamma}R(X)_\alpha
\]
where $R(X)_\alpha=\Gamma(X,M_\alpha)$ is the space of global sections of
$M_\alpha$.

In the case where $F = Cl(X)$ (so that $Cl(X)$ is free and finitely generated)
we write $R(X)$ for $R_F(X)$. The interesting issue here is to find useful conditions
on $X$ and $F$ which guarantee that $R_F(X)$ is finitely generated as a $K$-algebra.
In case $Cl(X)$ is finitely generated but not free, there is often an ``obvious"
free subgroup $F\subseteq Cl(X)$ of finite index that one would like to
consider in the discussion.

\begin{definition}
We say that the projective variety $X$ has a {\em factorial hull} if 
$Cl(X)$ is finitely generated and free, and $R(X)$ is finitely generated as a $K$-algebra. 
If $Cl(X)$ if merely finitely generated and there is a free abelian subgroup 
$F\subseteq Cl(X)$ of finite index such that $R_F(X)$ is a finitely generated 
$K$-algebra, we refer to $R_F(X)$ as a {\em $\mathbb{Q}$-factorial hull} of $X$.
\end{definition}

In this paper it will often suffice to think in terms of a $\mathbb{Q}$-factorial hull, since our 
main result is the calculation of a certain rational poyhedral cone $\mathscr{H}$ in the class group of $X$.
See Theorems~\ref{Hforthediva.thm} and ~\ref{fisacube.thm}. 

$\Gamma(X,M_\alpha)$ can often be described directly in terms of divisors
and functions on $X$. If $X$ has a factorial hull then $X$ can be described as the categorical
quotient of an open subset $U$ of $Spec(R(X))$ by an algebraic torus.
Furthermore, by the results of \cite{E-K-W}, $R(X)$ is a unique
factorization domain, hence our terminology {\em factorial hull}. Other
terminologies here are {\bf Cox ring} \cite{Cox}, and {\bf total homogeneous coordinate ring}
\cite{E-K-W}. Hausen \cite{Hausen2} has shown that, in some cases, this factorial hull 
can be used to provide greater flexibility in geometric invariant theory by allowing certain
Weil divisors instead of restricting only to ample Cartier divisors. $\mathbb{Q}$-factorial 
projective varieties with factorial hull are characterized in \cite{HK} using Mori theory.
It is known that any spherical variety has a $\mathbb{Q}$-factorial hull. However, nobody seems to have
published any proof. Brion has described a construction that depicts each $\mathbb{Q}$-factorial
hull as a fibre product over the appropriate versal object, similar to how Vinberg \cite{Vi} obtains
all flat, reductive monoids from $Env(G_0)$. Indeed, $Env(G_0)$ is a $\mathbb{Q}$-factorial 
hull of the {\bf plongement magnifique} \cite{DP}. Alternately, one can use
a result of Knop \cite{Knop} to obtain the $\mathbb{Q}$-factorial hull of a spherical
variety by observing that, for $X$ spherical, $R_F(X)$ is the ring of global functions on some
closely related spherical variety.

Our cone $\mathscr{H}$ appears to be what is called the {\bf moving cone} in \cite{HK}.
Unfortunately it is beyond the scope of this paper to discuss in detail the very interesting relationship
between our results and those of \cite{HK}. In fact our main results are independent of 
the discussion of factorial hulls or Mori theory. We mention those results mainly to place our work in a
more general context that could be appealing to some readers. Clearly there is potential for further
interesing work along these lines.

In some cases, (e.g. toric varieties \cite{Cox}, reductive monoids \cite{bible}, the $E_6$
cubic surface \cite{H-T}) $R(X)$ can be constructed geometrically from $X$ by realizing $Spec(R(X))$
as the total space of a type of universal ``$\mathbb{G}_m^r$-torsor" over $X$. In the case of 
reductive monoids a factorial hull had been already constructed indirectly in \cite{R2.6}.
See Theorem~\ref{classcover.thm} below.

Now let $M$ be a {\bf semisimple} monoid \cite{bible} with unit group $G$. By definition,
$M$ is reductive and normal, $M$ has a zero element $0\in M$, and dim$Z(G)=1$. Let
$X=X_M=(M\backslash \{0\})/K^*$ be the associated {\bf semisimple variety} \cite{R3}.
See also \cite{Tim}. 

We now describe the $\mathbb{Q}$-factorial hull of $X = X_M$. It suffices to descibe the
$\mathbb{Q}$-factorial hull of $M$ (whose coordinate algebra is a ``sufficiently good" summand of 
$R(X)$ so that $M$ and $X(M)$ have the same $\mathbb{Q}$-factorial
hull). For the details of this construction we refer to Theorem 6.7 of \cite{bible}.
We let $\Lambda^1(M)$ be the set of codimension-one $G\times G$-orbits of $M$.

\begin{theorem}                         \label{classcover.thm}
Let $M$ be reductive with unit group $G$. Then there exists a reductive monoid
$\widehat{M}$ and a morphism $\pi : \widehat{M}\to M$ such that
\begin{enumerate}
\item[(i)] $Cl(\widehat{M})=(0)$,
\item[(ii)] $\pi$ induces a bijection of\; $\Lambda^1
(\widehat{M})\to\Lambda^1(M)$.
\end{enumerate}
Furthermore, the unit group $\widehat{G}$ of $\widehat{M}$ is an extension of
$G$ by a $D$-group. 
\end{theorem}

\begin{remark}
$\widehat{M}$ was first constructed in \cite{R2.6} before the author noticed that 
such monoids are necessarily {see \bf flat} (\S~\ref{flat.sec}). The main problem
originally motivating the present paper is the following question.\\

\noindent {\em Let $M$ be a semisimple monoid with unit group $G$. How does one 
characterize/classify, up to isomorphism, the set of semisimple monoids $N$ such that 
$\widehat{N}=\widehat{M}$? What finer structure does this set $\mathscr{H}$ of monoids have?} \\

\noindent Although these questions might seem far removed from anything geometric, the final 
answer is most naturally phrased in terms of the divisor class group of
$M$ or $X_M=(M\backslash\{0\})/K^*$,
since this class group can be identified with a group of characters on the center $\widehat{M}$.
See Theorems~\ref{findh.thm}, ~\ref{Hforthediva.thm} and ~\ref{fisacube.thm}. Theorem~\ref{findh.thm}
discusses the case of a general semisimple monoid,
while Theorems~\ref{Hforthediva.thm} and ~\ref{fisacube.thm} discuss the case of the 
monoids (canonical monoids) associated with the plongement magnifique. 
\end{remark}

\begin{remark}     \label{ambighull.rk}
The construction of $\widehat{M}$ involves a (noncanonical) finite, dominant morphism
$\pi : G'\to G$ with $Cl(G')=(0)$. Since the latter is not canonical, neither is $\widehat{M}$.
To obtain a canonical object one needs to work a little deeper,
and also be satisfied with a $\mathbb{Q}$-factorial hull. 
\end{remark}

\begin{remark}         \label{qfact.ex}
Let $M$ be a semisimple monoid with associated semisimple variety
$X_M=(M\backslash\{0\})/K^*$. Let $\Lambda^1=\{[D_i]\}\subseteq Cl(X_M)$ be the
set of divisor classes of codimension-one $G\times G$-orbits of $X_M$.
Let $V\subseteq Cl(X_M)$ be the subgroup generated by
$\Lambda^1$. Then $V$ is a torsion-free subgroup of $Cl(X)$ of finite index.
Furthermore, $Cl(M)=V\oplus Cl(G)$.
$R_V(X_M)=R_V(M)$ is a $\mathbb{Q}$-factorial hull of $X_M$. 
Furthermore $K[M]\subseteq R_V(X_M)$ and $M_V=Spec(R_V(X_M))$ is a 
reductive normal algebraic monoid whose unit group is an extension
of the unit group of $M$ by an algebraic torus.
See \S 2 below for another construction, called $M_\mathbb{H}$, in 
terms of the hyperplanes in $X(T_0)$ determined by 
$\{[D_i]\}$. It turns out that there is a canonical isomorphism 
$M_V\to M_\mathbb{H}$.
\end{remark}

\section{Flat Monoids and Hyperplane Arrangements}  \label{flat.sec}

Let $M$ be a semisimple monoid. In this section we discuss the construction 
$M\leadsto M_\mathbb{H}$ (and related constructions) from several points of view. 
On the one hand, it turns out that $M_\mathbb{H}$ is determined by a certain rational, 
oriented, hyperplane arrangement in $X(T_0)$ (a {\bf $W$-arrangment}). On the 
other hand, any reductive monoid with trivial divisor class group is 
{\bf flat} in the sense of Vinberg \cite{Vi}. We determine how the 
$W$-arrangement $\mathbb{H}$ of $M$ essentially determines the monoid 
$M_\mathbb{H}$. Indeed, the theory of flat monoids allows us to construct 
$M_\mathbb{H}$ directly from $\mathbb{H}$.

\subsection{Flat Reductive Monoids}                 \label{flat.sub}

Associated with any reductive monoid $M$, is its {\bf abelization}
\[    
\pi : M\to A.    
\]
In \cite{Vi} Vinberg calls $M$ {\bf flat} if $\pi$ is a flat morphism with
reduced and irreducible fibres. An important observation here 
(see Theorem~\ref{vinbergflat.thm} below) is that any reductive monoid 
$M$ with trivial divisor class group is flat.

Let $G$ be the unit group of $M$, and let $G_0$ be the
semisimple part of $G$. Let $B$ and $B^-$ be opposite Borel subgroups
of $G$ containing the maximal torus $T$ with unipotent radicals
$B_u$ and $B_u^-$ respectively. Let $T_0=T\cap G$ be the 
associated maximal torus of $G_0$. 

Let $X(T_0)_+$ denote the monoid of dominant weights of $T_0$. If
$\lambda\in X(T_0)_+$, we can write

\[    \lambda = \sum_{\alpha\in\Delta}c_{\alpha}\lambda_{\alpha},     \]
where $\{\lambda_{\alpha}\}$ is the set of fundamental dominant weights
of $G_0$. Define

\[         c : X(T_0)_+\to Cl(M)         \]
by
$c(\lambda)=\sum_{\alpha\in\Delta}c_\alpha[\overline{Bs_{\alpha}B^-}]$.
Here $s_\alpha\in S$ is the simple involution corresponding to 
$\alpha\in\Delta$. Let
\[
L(M)=\{f\in K[M]\;|\; f(ugv)=f(g)\;\text{for all}\; u\in B_u,\;v\in B_u^-
\;\text{and}\;f(1)=1\;\}.
\]
Let $Z\subseteq G$ be the connected center of $G$ so that $G=ZG_0$ and let 
$\overline{Z}\subseteq M$ be the Zariski closure of $Z$ in $M$. $X(\overline{Z})$
is the set of characters of $\overline{Z}$.

\begin{theorem}                     \label{vinbergflat.thm}
Let $M$ be a reductive monoid with unit group $G$, and let $G_0$ be the
semisimple part of $G$. Assume that $M$ has a zero element. 
The following are equivalent.
\begin{enumerate}
\item[a)] The abelization morphism $\pi : M\to A$ is flat, with reduced and
irreducible fibres. 
\item[b)] The following two conditions hold.
  \begin{enumerate}
   \item[i)] If $\chi_1\lambda_1=\chi_2\lambda_2$ ($\lambda_i\in\mathscr{M}$,
   $\chi_i\in X(\overline{Z})$) then $\chi_1=\chi_2$ and $\lambda_1=\lambda_2$.
   \item[ii)] $\mathscr{M}$ is a subsemigroup of $L(M)$. 
  \end{enumerate}
\item[c)] The canonical map $c : X(T_0)_+\to Cl(M)$ is trivial. 
\item[d)] For any irreducible representation $\rho : M\to End(V)$ there
is a character $\chi : \overline{Z}\to K$ of $\overline{Z}$, and an irreducible
representation $\sigma : M\to End(V)$, such that $\sigma(e)\neq 0$ for any $e\in
\Lambda^1$ and $\rho=\chi\otimes\sigma.$   
\item[e)] Any $f\in L(M)$ factors as $f=\chi g$ where $\chi\in X(\overline{Z})$, and
$g\in\L(M)$ has zero set
$Z(g)\subseteq \cup_{\alpha\in\Delta}\overline{Bs_{\alpha}B^-}$.
\end{enumerate}
\end{theorem}

For the proof see Theorem 6.10 of \cite{bible}. Notice in particular, 
if the divisor class group of $M$ is trivial, that $M$ is flat (using part c)
of Theorem~\ref{vinbergflat.thm} above). 

It turns out that that there is a universal, flat monoid $Env(G_0)$ associated
with each semisimple group $G_0$. This amazing monoid was originally discovered
and constructed by Vinberg in \cite{Vi}. He refers to it as the {\bf
enveloping semigroup} of $G_0$. It has the following universal property.

Let $M$ be any flat monoid with zero. Assume that the
semisimple part of the unit group of $M$ is $G_0$. Let $A(M)$ denote the
abelization of $M$, and let $\pi_M : M\to A(M)$ be the abelization morphism.
We make one exception with this notation. We let $A$ denote the abelization
of $Env(G_0)$ and we let $\pi : Env(G_0)\to A$ be the abelization morphism.
The universal property of $Env(G_0)$ is as follows. Given any isomorphism
$\varphi_0$ from the semisimple part of $G(M)$ to the semisimple
part of $Env(G_0)$, there are unique morphisms

\[        a : A(M)\to A                    \]
and

\[        \varphi : M\to Env(G_0)          \]
such that 
\begin{enumerate}
\item[i)]  $\varphi|G_0=\varphi_0$;
\item[ii)] $a\circ\pi_M = \pi\circ\varphi$; 
\item[iii)] $\phi : M\cong E(a,\pi)$, via $\phi(x)=(\pi_M(x),\varphi(x))$, where
$E(a,\pi)=\{(x,y)\in A(M)\times Env(G_0)\;|\;a(x)=\pi(y)\}$, is the 
{\bf fibred product} of $A(M)$ and $Env(G_0)$ over $A$. 
\end{enumerate}

There are several ways to construct this monoid $Env(G_0)$, and there
are already hints in Theorem~\ref{vinbergflat.thm}. However, we use
the construction in Theorem 17 of Rittatore's thesis \cite{Rit0}. 
The reader should also see Vinberg's construction in \cite{Vi}. A
similar construction, due to Rittatore \cite{Rit1}, exists in positive 
characteristics. Notice that we are using multiplicative notation for 
characters. In particular, $X(A)\cong P$ is the submonoid of $X(T_0)$
generated by the positive roots. 

Assume that $M$ has unit group $G$. Since $G\subseteq M$ 
is open we obtain that
\[  K[M]\subseteq K[G].  \]
Now, it is well known that
\[ K[G] = \bigoplus_{\lambda\in X_+} K[G]_\lambda \]
where $X_+$ is the set of dominant characters of $T$ with respect to
$B$.  Here, each $K[G]_\lambda$ is an irreducible $G\times G$-module 
with highest weight $\lambda\otimes\overline{\lambda}$ and  $G\times G$
acts on $G$ via $((g,h),x)\mapsto gxh^{-1}$.  This ``multiplicity
$\leq 1$" condition implies that any $G\times G$-stable subspace of
$K[G]$ is a sum of some of the $K[G]_\lambda$.  In particular,
\[
K[M] = \bigoplus_{\lambda\in L(M)} K[G]_\lambda,
\]
where $L(M)\subseteq X_+$.  We refer to $L(M)$ as the {\bf augmented
cone} of $M$. It is defined, in a different but equivalent way, in 
Section~\ref{flat.sub}.
 
\begin{theorem}                    \label{envg0.thm}
Let $G_0$ be a semisimple group and let
\[
\mathscr{L}(G_0)=\{(\chi,\lambda)\in L(T_0\times G_0)\;|\;
\chi\lambda^{-1}\in P\}.
\]
Define
\[
K[Env(G_0)]=\bigoplus_{(\chi,\lambda)\in \mathscr{L}(G_0)}(V_\lambda\otimes
V_\lambda^*)\otimes\chi\;\subseteq K[G_0\times T_0].
\]
\noindent Then $K[Env(G_0)]$ is the coordinate algebra of the normal, 
reductive algebraic monoid $Env(G_0)$ with the above-mentioned universal 
property. In particular, $L(Env(G_0))=\mathscr{L}(G_0)$.
\end{theorem}

For the proof see Theorem 6.16 of \cite{bible}.

If $G_0$ is simply connected then $K[Env(G_0)]$ is the factorial hull
of any {\bf canonical monoid} of $G_0$ \cite{PR}.

Closer scrutiny of this universal property yields a 
numerical classification of the flat monoid $M$ in terms of a certain 
map $\theta_M^* : X(T_0)_+\to X(\overline{Z})$. The above-mentioned 
classifying map  $a : A(M)\to A$ can then be calculated directly.

We view $X(T_0)$ as a multiplicative group and write the simple roots
exponentially $\{e^\alpha\;|\;\alpha\in\Delta\}\subseteq X(T_0)$.
Recall that $X(A)\cong P$, the free commutative submonoid
of $X(T_0/Z_0)$ generated by the positive roots. 
Thus the coordinate ring of $A$ is the polynomial ring with the universal 
generators $\{u_\alpha\cong e^\alpha\;|\;\alpha\in\Delta\}$.

\begin{corollary}     \label{classifmap.cor}
Let $M$ be flat with unit group $G=G_0Z$ and connected center 
$Z\subseteq G$. Let $Z_0=Z\cap G_0$.
\begin{enumerate}
 \item $M$ is determined by a certain map $\theta_M^* : X(T_0)_+\to X(\overline{Z})$
 such that
  \begin{enumerate}
   \item $\theta|Z_0=id$,
   \item $\theta_M^*$ extends to $\theta_M^* : X(T_0)\to X(Z)$ with 
   $\theta_M^*(e^\alpha)\in X(A_M)=X(\overline{Z}/Z_0)\subseteq X(\overline{Z})$ for all 
   $\alpha\in\Delta$.
  \end{enumerate}
  Conversely, any $\theta^*: X(T_0)_+\to X(\overline{Z})$ satisfying 
  a) and b) above determines a flat monoid
  $M=M_\theta$ with unit group $G$.
 \item The augmented cone of $M_\theta$ is determined by $\theta$ as follows.
 \[
 L(M_\theta)=\{(\chi,\lambda)\in L(Z\times G_0)\;|\;\chi\theta^*_M(\lambda)^{-1}\in X(A)\}.
 \]
 \item To obtain $M$ as a fibred product from $\pi : Env(G_0)\to A$ define
\[
a : A_M\to A
\]
by the rule 
\[
a^*(u_\alpha)=\theta_M^*(e^\alpha).
\]
Then, as above, $\phi : M_\theta\cong E(a,\pi)$ via $\phi(x)=(\pi_M(x),\varphi(x))$.
\end{enumerate}
\end{corollary}

The above corollary is a reformulation of Theorems 4 and 5 of \cite{Vi}.
See also the proof of Theorem 6.16 in \cite{bible}. Recall from Theorem~\ref{envg0.thm}
that 
\[
L(Env(G_0))=\mathscr{L}(G_0)=\{(\chi,\lambda)\in L(T_0\times G_0)\;|\;
\chi\lambda^{-1}\in P\}.
\]
In this case, $X(\overline{Z})=\{\lambda\in X(T_0)\;|\;\lambda^n\in P
\;\text{for some}\; n>0\}$ and $\theta^* : X(T_0)_+\to X(Z)$ is just the inclusion.

There is another useful characterization of flat monoids. 
Let $M$ be a reductive monoid and let $\overline{T}\subseteq M$ be the closure 
in $M$ of a maximal torus of $G$. Define 

\[
X(\overline{T})_+=\{\chi\in X(\overline{T})\;|\;\Delta_\alpha(\chi)\geq0\;\text{for all}\;\alpha\in\Delta\},
\]
where $\Delta_\alpha$ is defined by the equation 
\[
\Delta_\alpha(\chi)\alpha=\chi - \sigma_\alpha(\chi).
\]
For each $e\in\Lambda^1$ there is a ``valuation" $\nu_e : X(T)\to\mathbb{Z}$
determined by the divisor $e\overline{T}\subseteq\overline{T}$ ($\nu_e$ is 
induced by the inclusion $K^*\subseteq T$ of the 1-PSG containing $e$ in its closure).

Define
\[
\mathscr{M}=\{\lambda\in X(\overline{T})_+\;|\;\nu_e(\lambda)=0\;\text{for all}\; e\in\Lambda^1\}.
\]

\begin{theorem} \label{moreflat.thm}
Let 
\[
r : \mathscr{M}\to X(T_0)_+
\]
be defined by $r(\lambda)=\lambda|T_0$.
The following are equivalent.
\begin{enumerate}
\item $r : \mathscr{M}\to X(T_0)_+$ is an isomorphism.
\item $M$ is flat
\end{enumerate}
Furthermore, in this case, $\theta^*=p\circ r^{-1}$, where 
$p : \mathscr{M}\to X(\overline{Z})$ is defined by $p(\lambda)=\lambda|\overline{Z}$.
\end{theorem}

\begin{proof}
Assume that $M$ satisfies condition {\em 1} above.
If $\lambda\in X(T)_+$ then there is a unique
$\lambda_0\in\mathscr{M}$ such that $r(\lambda)=r(\lambda_0)$. So we let
$\delta = \lambda\lambda_0^{-1}$. It follows easily that
$\delta\in X(A)\subseteq X(\overline{T})$. This gives us the
desired factorization $\lambda=\delta\lambda_0$ as in part e) of
Theorem~\ref{vinbergflat.thm}.
\end{proof}

\subsection{The Flat Monoid $M_\mathbb{H}$ of an Arrangement} 

Let $M$ be semisimple monoid with unit group $G$ and let $T$ be a maximal 
torus of $G$. Recall that
$E(\overline{T})=\{e\in\overline{T}\;|\;e^2=e\}$. $E(\overline{T})$ is a
poset if we define $e\geq f$ whenever $ef=f$. Let
\[
E^1=\{e\in E(\overline{T})\;|\;e\;\text{is maximal in}\;E(\overline{T})\backslash\{1\}\}.
\]
$E^1$ is refered to as the set of {\bf maximal idempotents} of $\overline{T}$.
$E^1$ corresponds bijectively to the set of codimension-one 
$T\times T$-orbits $\{eT\}$ of $\overline{T}$.
Associated with each $eT\subseteq\overline{T}$ 
there is the associated valuation $v_e$ of $K[T]$.
Let $Z\subseteq G$ be the connected center of $G$. In the following
definition we let $\nu_e$ be the restriction of $v_e$ to $X(T/Z)\subseteq X(T)$.
$\nu_e$ is dual to the 1-psg $\lambda : K^*\to T/Z$ that has 
$\underset{t\to0}{lim}\lambda_t\in (eT)/Z$.

\begin{definition}  \label{arrange.def}
Let $M$ be as above.
\begin{enumerate}
\item [a)]
The {\em arrangement} $\mathbb{H}(M)$ of $M$ is the collection 
\[
\mathbb{H}(M) = \{\nu_e\;|\; e\in E^1\}.
\]
\item [b)] Let
\[
\Lambda^1\mathbb{H}(M)=\{\nu\in\mathbb{H}(M)\;|\;\nu(-\alpha)\geq 0
\;\text{for all}\;\alpha\in\Delta\}.
\]
\end{enumerate}
\end{definition}

\noindent $\Lambda^1\mathbb{H}(M)$ is a fundamental domain for the action of the
Weyl group on $\mathbb{H}(M)$. $\Lambda^1\mathbb{H}(M)$ can be identified with 
$\Lambda^1=\{e\in E^1(\overline{T})\;|\;Be=eBe\;\}$.

We identify $\mathbb{H}(M)$ and 
$\mathbb{H}(N)$ if there is a bijection $f : \mathbb{H}(M)\to\mathbb{H}(N)$ such that, 
if $f(l)=m$, then $ker(l)=ker(m)$ and $l$ is a positive multiple of $m$. 
Notice also that if $l,m\in\mathbb{H}(M)$ and $ker(l)=ker(m)$, and $l$ is a 
{\em positive} multiple of $m$, then $l=m$. This is a natural 
nondegeneracy condition inherent in a semisimple monoid. We conclude that 
$\mathbb{H}(M)$ can be thought of as a collection of rational, oriented, 
$W$-invariant hyperplanes $\{(H,l)\}$ in $X(T/Z)\otimes\mathbb{Q}$. 

By abuse of language we often write $\mathbb{H}(M) = \{(ker(\nu_e),\nu_e)\;|\; e\in E^1\}$,
if we wish to emphasize the r\^ole of $ker(\nu_e)$. In any case, 
$\mathbb{H}(M)$ can be thought of as either a set oriented hyperplanes in $X(T/Z)$, 
or as a set of rays in $Hom_\mathbb{Z}(X(T/Z),\mathbb{Z})$.

If $w\in W$ the action on $\mathbb{H}(M)$ is determined by

\[
w(H,l)=(w(H),l\circ w^{-1}).
\]
Notice also that we may identify $X(T/Z)\otimes\mathbb{Q}$ with 
$X(T_0)\otimes\mathbb{Q}$, where $T_0\cap G_0\subseteq G_0$ is the 
associated maximal torus of the the semisimple part of $G$.

We now construct $M_\mathbb{H}$ from an arrangement $\mathbb{H}$ using 
Corollary~\ref{classifmap.cor}. Let $\Lambda^1\mathbb{H}$ be 
a finite set of hyperplanes in $X=X(T_0)$ such that for each 
$H\in\Lambda^1\mathbb{H}$, there is a functional $l : X\to\mathbb{Z}$ such that 

\begin{enumerate}
\item $l(-\alpha)\geq 0$ for all $\alpha\in\Delta$,
\item $ker(l)=H$.
\end{enumerate}

Notice that $Span_{\mathbb{Z}}(\Delta)=X(T/Z_0)\subseteq X$ is a subgroup of finite 
index. Notice also that $l$ is determined up to a positive scalar by $H$. So 
each $H$ is canonically oriented. We then define $\mathbb{H}$ as follows.
\[
\mathbb{H}=\{(w(H),l\circ w^{-1})\;|\;w\in W\}.
\]
Now for each $H\in\Lambda^1\mathbb{H}$, there is a unique (``best") functional
$l : X\to\mathbb{Z}$, as above, which is not an integer multiple of any 
other. So we can write $(H,l)$ for $H$. Write $\Lambda^1\mathbb{H}=\{(H_1,l_1),(H_2,l_2),...,(H_s,l_s)\}$. 

Let $\pi : Env(G_0)\to A$ be the abelization morphism of $Env(G_0)$.
Notice that the set of characters $X(A)$ of $A$ is canonically identified with
the submonoid $<\Delta>$ of $X(A^*)$. Thus, 
$<\Delta>=P=X(A)\subseteq X(T_0/Z_0)=X(A^*)=Span_\mathbb{Z}(\Delta)\subseteq X(T_0)=X$.
Hence each functional $l : X\to\mathbb{Z}$ restricts to a functional $l : X(A)\to\mathbb{Z}$.
Define 
\[
a^* : X(A)\to\mathbb{N}^s  
\]
by $a^*(\chi)=(-l_1(\chi),-l_2(\chi),...,-l_s(\chi))$. This makes sense, since 
for each $\alpha\in\Delta$ and each $i$, $l_i(-\alpha)\geq 0$.
Define
\[
X(A_{\mathbb{H}})=\oplus_{i=1}^sk_i\mathbb{N}\subseteq\mathbb{N}^s
\]
where each integer $k_i>0$ is chosen to be maximal subject to the condition
that $a^*(X(A))\subseteq X(A_{\mathbb{H}})$.

Thus we obtain, by definition, that
\[
a : A_{\mathbb{H}}\to A.
\]

\begin{definition}  \label{arrangemonoid.def}
We let $M_{\mathbb{H}}=E(a,\pi)$ be the fibre product  
of $a$ over $\pi$ (as in Corollary~\ref{classifmap.cor}). 
$M_{\mathbb{H}}$ is called the {\em flat monoid of the 
arrangement} $\mathbb{H}$.
\end{definition}

\begin{remark}  \label{thetastarvsastar.rk}
The purpose of this remark is to indicate in more detail how $\theta^*$ 
and $a^*$ determine each other. Suppose that $M$ is flat. 
Then for all $e\in E^1(A)$ there is a unique $\chi\in X(A_M)$
such that $eA_M=\chi^{-1}(0)$. Furthermore, this defines a ``valuation"
$\nu_e :L(G)\to \mathbb{Z}$ since $\pi^{-1}(eA_M)$ is a codimension-one
$G\times G$-orbit of $M$. Now
\[ 
\mathscr{M}=\{(\theta^*(\lambda),\lambda)\;|\;\lambda\in X(T_0)_+\}.
\]
Thus,
\[
\nu_e((1,\lambda))=\nu_e((\theta^*(\lambda)^{-1},1)(\theta^*(\lambda),\lambda))
=\nu_e(\theta^*(\lambda)^{-1},1)
\]
since, by Theorem~\ref{moreflat.thm}, $\nu_e(g)=0$ for all $g\in\mathscr{M}$ and 
$e\in E^1(A_M)$. Observe also that 
\[
\nu_e((1,e^{-\alpha}))=\nu(\theta^*(e^{\alpha}),1)\geq 0
\]
since $(\theta^*(e^{\alpha}),1)\in X(A_M)\subseteq X(\overline{Z})$.

If $M\cong E(a,\pi)$, as in Definition~\ref{arrangemonoid.def}, one can check that 
$a^* : P\to X(A_M)$ extends uniquely to a map $b^* : X(Z(Env(G_0)))\to X(\overline{Z})$.
Furthermore,
\[
L(M)=\{(\delta b^*(\lambda),\lambda)\;|\;\lambda\in X(T_0)_+, \delta\in X(A_M)\}.
\]
The canonical map $M\cong E(a,\pi)\to Env(G_0)$ induces the map 
$\gamma : L(Env(G_0))\to L(M)$ defined by $\gamma(\chi,\lambda)=(b^*(\chi),\lambda)$.
\end{remark}

Recall from  Example~\ref{qfact.ex} the definition of $M_V$. 

\begin{theorem}  \label{mvisisotomh.thm}
There is a canonical finite, dominant morphism $M_V\to M_\mathbb{H}$.
\end{theorem}
\begin{proof}
The canonical restriction map $Cl(M_V)\to Cl(G)$ is an isomorphism, 
by construction of $M_V$, since $ker(Cl(M_V)\to Cl(G))$ is the 
subgroup generated by $\{[D_i]\}$ as in Remark~\ref{qfact.ex}. 
However, by construction of $M_V$, each $[D_i]\in Cl(M_V)$ is zero.
Also the composite $X(T_0)_+\underset{c}{\to} Cl(M)\to Cl(G)$ is 
the trivial homomorphism for any reductive, normal monoid. Thus
$c : X(T_0)_+\to Cl(M_V)$ is also trivial. So $M_V$ is flat by 
Theorem~\ref{vinbergflat.thm}. Thus by Corollary~\ref{classifmap.cor}
there is a unique classifying map
\[
h^* : X(A)\to X(A_{M_V})
\]
such that $M_V\cong E(h,\pi)$. The assumption that $Cl(M_V)=Cl(G)$ is
equivalent to the statement that for each $e\in\Lambda^1$ the ideal
$\{f\in K[M_V]\;|\; f\;\text{vanishes on}\; GeG\}$ is a principal ideal
$(\chi_e)$ of $K[M]$. Hence $X(A_M)=\mathbb{N}^s$ is the free abelian
monoid on these (normalized) generators. Let $p_e : X(A_M)\to\mathbb{N}$
be the projection onto the $i$-th factor and let $h_e^*=p_e\circ h^*$.
Let $X\subseteq X(A_{M_V})$ be the smallest submonoid of $X(A_M)$
of the form
\[
X=\oplus_{e\in\Lambda^1}k_e\mathbb{N}
\]
so that $h^*(X(A))\subseteq X$. Since $M_V$ is a monoid with 
$|\Lambda^1|=|\mathbb{H}|$ one can check that
$h^* : X(A)\to X$ is isomorphic to $a^* : X(A)\to X(A_\mathbb{H})$,
and that $K[X]\subseteq K[X(A_M)]$ is a finite morphism. 
This says that the classiying map for $M$ factors through the classifying
map for $M_\mathbb{H}$. Thus there is a canonical finite morphism
$M\cong E(h,\pi)\to M_\mathbb{H}$. Furthermore $E(h,\pi)\cong E(a,\pi_\mathbb{H})$.
\end{proof}

\begin{corollary}    \label{mvisisotomh.cor}
Let $M$ be a reductive, normal monoid with unit group $G$ and zero element
$0\in M$. Suppose that $Cl(M)=Cl(G)$. Then $M$ is flat.
\end{corollary}

\section{The Coterie of an Arrangement}   \label{amp.sec}

In this section we shall often work with rational cones and vector spaces. 
But we want our notation to be consistent with all previous notation. If $L$
is some lattice of interest we denote by $L_0$ the rational vector space 
$L\otimes_{\mathbb{Z}}\mathbb{Q}$. If $P\subseteq L$ is a finitely 
generated, additive submonoid of $L$ then $P_0$ is the rational cone
in $L_0$ generated by $P$ (In this paper $P_0$ is the rational polyhedral cone 
generated by the positive roots). We make one exception with this notation. 
Script quantities (esp. $\mathscr{H}$, $\mathscr{C}$ and $\mathscr{C}(\delta)$;
see below) automatically represent rational entities. 

The purpose of this section is to determine a useful criterion
for identifying a certain rational polyhedral cone $\mathscr{H}$. 
This cone $\mathscr{H}$ can be thought of as the
Weil divisor analogue of the ample cone. In the case of general $\mathbb{H}$ we
obtain, in a straight forward manner, that $\mathscr{H}\subseteq
X(\overline{Z(\mathbb{H})})_0$.

\subsection{The Cone $\mathscr{H}$} \label{H.subsec}

Throughout this section we use additive notation for any calculations with
characters. Scalars are less cumbersome than exponents 
since we are are working with rational vector spaces.

Let $M$ be a semisimple monoid with unit group $G$ and semisimple part
$G_0\subseteq G$. By the results of \S~\ref{flat.sec} 
we have a canonical quotient morphism
\[
M_V\to M
\]
and a canonical finite morphism
\[
M_V\to M_\mathbb{H}
\]
where $\mathbb{H}=\mathbb{H}(M)$. Since we shall be working in the appropriate 
rational cone, there is no harm in assuming that we actually have a morphism
\[
M_\mathbb{H}\to M.
\]
This amounts to replacing the $\delta$ in Proposition~\ref{lofmdelta.prop} below by 
$k\delta$ for some $k\in\mathbb{N}$. In any case, we can restrict characters, functions 
etc. on $M$ to obtain like quantities on $M_\mathbb{H}$. For example, if $\delta\in X(M)$, 
we can write $\delta\in X(M_\mathbb{H})$.

Let $T_0\subseteq G_0$ be a maximal torus and let $\overline{Z(\mathbb{H})}\subseteq
M_\mathbb{H}$ be the Zariski closure of the connected center $Z(\mathbb{H})$ of the
unit group of $M_\mathbb{H}$. Recall also the abelization $\pi : M_\mathbb{H}\to A_\mathbb{H}$ 
of $M_\mathbb{H}$. We can identify $X(A_\mathbb{H})$ as a submonoid of
$X(\overline{Z(\mathbb{H})})$ by restriction of $\pi$. By Corollary~\ref{classifmap.cor} there
is a structure map for $M_\mathbb{H}$,

\[
\theta^* : X(T_0)\to X(\overline{Z(\mathbb{H})}).
\]
Recall from the proof of Theorem 6.16 of \cite{bible},
that 
\[
L(M_\mathbb{H})=\{(\gamma,\lambda)\in X(\overline{Z})\oplus
X(T_0)\;|\; \gamma-\theta^*(\lambda)\in X(A_\mathbb{H})\}.
\]

\begin{proposition}   \label{lofmdelta.prop}
There is a unique $\delta\in X(\overline{Z(\mathbb{H})})$ such that
\[
L(M)=\{(k\delta,\lambda)\in X(\overline{Z(\mathbb{H})})\oplus
X(T_0)\;|\;k\delta-\theta^*(\lambda)\in X(A)\;\text{and}\; k\geq
0\}.
\]
\end{proposition}
\begin{proof}
It is straightforward to check that 
\[
L(M)=L(M_\mathbb{H})\cap \{(\gamma,\lambda)\in X(\overline{Z(\mathbb{H})})\oplus X(T_0)\;|\;
\gamma=k\delta\;\text{for some $k\geq 0$}\;\}. 
\]
$\delta$ is the generator of $X(M)\cong\mathbb{N}$, the character monoid of $M$.
\end{proof}

It is convenient to work over the positive rational numbers. Define

\[
\mathscr{C}=X(T_0)_+\otimes\mathbb{Q}^+,
\]
the {\em rational Weyl chamber}.

\begin{definition}  \label{crosssectionpolytope.def}
If $\delta\in X(\overline{Z(\mathbb{H})})$ is as in Proposition~\ref{lofmdelta.prop}
above we define
\[
\mathscr{C}(\delta)=\{\lambda\in\mathscr{C}\;|\;\delta-\theta^*(\lambda)\in
X(A)_0\}.
\]
$\mathscr{C}(\delta)$ is called the {\em cross section polytope of $M$}.
\end{definition}

$\mathscr{C}(\delta)$ is a fundamental domain for the action of the Weyl
group $W$ of $G_0$ on the polytope 

\[
\mathscr{P}(\delta)=\cup_{w\in W}w(\mathscr{C}(\delta))\subseteq X(T_0)\otimes\mathbb{Q}.
\]
As above let $\mathbb{H}=\mathbb{H}(M)$ be the arrangement of $M$ and let $M_\mathbb{H}$ be 
the associated flat monoid with classifying map
\[
\varphi : M_\mathbb{H}\to M.
\]
Since $\varphi$ is surjective $\varphi^* : K[M]\to K[M_\mathbb{H}]$ identifies
$K[M]$ as a subalgebra of $K[M_\mathbb{H}]$. The action of $Z(\mathbb{H})$ on 
$M_\mathbb{H}$ determines a direct sum decomposition 
\[
K[M_\mathbb{H}]=\oplus_{\chi\in X(\overline{Z(\mathbb{H})})}K[M_\mathbb{H}]_\chi.
\]
Let $Z\subseteq G$ be the connected center of the unit group of $M$. 
Since $dim(Z)=1$, $K[M]$ is graded as follows.
\[
K[M]=\oplus_{n\geq 0}K[M_\mathbb{H}]_{n\delta}
\]
for some unique $\delta\in X(Z)\subseteq X(Z(\mathbb{H}))$. $\delta$ is the 
generator of $X(\overline{Z})\cong\mathbb{N}$.

The major problem here is to characterize the subset of $X(\overline{Z(\mathbb{H})})$
consisting of all possible characters $\delta$ that can arise from a semisimple monoid
$N$ with $\mathbb{H}(N)=\mathbb{H}$. While this might seem indirect and removed from 
the underlying geometry, it is {\em the} question that motivated the current work.
Proposition~\ref{gettingcofdelta.prop} and Theorem~\ref{findh.thm}
below allow us to calculate this interesting cone for flat monoids of the
form $M_\mathbb{H}$.

We denote the quadratic form on $\mathscr{C}$ by $(-,-)$.

\begin{proposition} \label{gettingcofdelta.prop}
For $\lambda\in\mathscr{C}$ we let $u^\delta_i(\lambda)=\nu_i(\delta-\theta^*(\lambda))$, where
$\{\nu_i : X(Z(\mathbb{H}))\to\mathbb{Z}\;|\; i\in I\}$ is the set of essential valuations
of $X(Z(\mathbb{H}))$ so that $X(\overline{Z(\mathbb{H})})=\cap_i\{\chi\;|\;\nu_i(\chi)\geq 0\}$. Then
\begin{enumerate}
\item
$\mathscr{C}(\delta)=\{\lambda\in\mathscr{C}\;|\;u^\delta_i(\lambda)\geq
0\;\text{for all $i$}\}$.
\item For each $i$, there is a unique $r_i(\delta)\in\mathscr{C}$ such that
$u^\delta_i(\lambda)=\epsilon_i(\delta)(r_i(\delta)-\lambda,r_i(\delta))$.
Furthermore,
$\epsilon_i(\delta)=\frac{\nu_i(\delta)}{(r_i(\delta),r_i(\delta))} >
0$.
\end{enumerate}
\end{proposition}
\begin{proof}
From Definition~\ref{crosssectionpolytope.def} we have that
\[
\mathscr{C}(\delta)=\{\lambda\in\mathscr{C}\;|\;\delta-\theta^*(\lambda)\in
X(A)_0\}.
\]
But also $X(A)_0=X(\overline{Z})_0=\{\chi\in X(Z)_0\;|\;\nu_i(\chi)\geq 0\}$.
Thus, 
\[
\mathscr{C}(\delta)=\{\lambda\in\mathscr{C}\;|\;u^\delta_i(\lambda)\geq
0\;\text{for all $i$}\}.
\]
This proves a).

Now let $H_i=\{\mu\in X(T_0)\;|\;\nu_i(\theta^*(\mu))=0\}$. Then 
$(H_i,l_i)\in\mathbb{H}$, where $l_i(\mu)=-\nu_i(\theta^*(\mu))$.
So let 
\[
K(u_i^\delta)=\{\lambda\in\mathscr{C}\;|\;u^\delta_i(\lambda)=0\}.
\]
Then there exists a unique $r_i(\delta)\in\mathscr{C}$ such that 
\begin{enumerate}
\item[a)] $\nu_i(\theta^*(r_i(\delta)))=\nu_i(\delta)$,
\item[b)] $(r_i(\delta),\mu)=0$ for all $\mu\in H_i$.
\end{enumerate}

We claim that $K(u_i^\delta)=r_i(\delta)+H_i$. Indeed, if $\mu\in H$,
then we get $u_i^\delta(r_i(\delta)+\mu)=0$ (using a) above and the 
fact that $\mu\in H_i$). Thus, $K(u_i^\delta)=r_i(\delta)+H_i$, since 
$K(u_i^\delta)\supseteq r_i(\delta)+H_i$, while they are both 
affine subspaces of the same dimension.

We claim also that there exists $\epsilon_i\in\mathbb{Q}^+$ such that 
$u^\delta_i(\lambda)=\epsilon_i(\delta)(r_i(\delta)-\lambda,r_i(\delta))$.
Indeed, if $\lambda=r_i(\delta)+\mu\in K(u_i^\delta)$ then 
\[
(r_i(\delta)-\lambda,r_i(\delta))=(-\mu,r_i(\delta))=0.\;\;\;\;\; (*)
\]
Then we let 
\[
\epsilon_i(\delta)=\frac{\nu_i(\delta)}{(r_i(\delta),r_i(\delta))}.
\]
Then for $\lambda\in r_i(\delta)+H_i$, 
\[
u^\delta_i(\lambda)=\epsilon_i(\delta)(r_i(\delta)-\lambda,r_i(\delta))
\]
since by (*) the RHS is zero, and by definition the LHS is zero.
But also, 
\[
u_i^\delta(0)=\nu_i(\delta)=\epsilon_i(\delta)(r_i(\delta)-0,r_i(\delta)).
\]
Thus, the formula is true for all $\lambda$.
\end{proof}

The following Theorem characterizes membership in the cone $\mathscr{H}$. It may not look like 
a ``Theorem" but it yields an exact method for calculating $\mathscr{H}$ from $\mathbb{H}$.

\begin{theorem}       \label{findh.thm}
Let $M=M_\mathbb{H}$. Then
$\mathscr{H}=\{\delta\in X(\overline{Z(\mathbb{H})})_0\;|\; \{r_i(\delta)\}\;\text{satisfies
the following condition}\}$.\\

For each $i$ there exists $x_i\in\mathscr{C}^0$, the interior of $\mathscr{C}$, such that 
\begin{enumerate}
\item $(r_j(\delta)-x_i,r_j(\delta)) > 0$ for all $j\neq i$.
\item $(r_i(\delta)-x_i,r_i(\delta)) = 0$.
\end{enumerate}
\end{theorem}
\begin{proof}
This is an exact reformulation, in terms of the $r_i$'s, of the 
(nondegeneracy) condition that assures us that none of the 
$H_i$'s is ``lost" when restricting from $M_{\mathbb{H}}$ to
$M$. It is equivalent to saying that, for each $i$, there is an
$x_i\in (r_i(\delta)+H_i)\cap\mathscr{C}^0$, that is ``below" the hyperplane
$r_j(\delta)+H_j$ for each $j\neq i$.
\end{proof}

\begin{corollary}  \label{hismultipl.cor}
$\mathscr{H}+\mathscr{H}\subseteq\mathscr{H}$.
Furthermore, for each $i$, $r_i(\delta+\gamma)=r_i(\delta)+r_i(\gamma)$.
\end{corollary}
\begin{proof}
The basic idea here is this. If ${\bf x}=(x_i)$ works for $\delta$, and ${\bf y}=(y_i)$
works for $\gamma$ then ${\bf x+y}$ works for $\delta+\gamma$.

Since both $\delta$ and $\gamma$ come from the same collection $\mathbb{H}$ we obtain that,
for each $i$,
\[
r_i(\gamma)=\alpha_ir_i(\delta)
\]
for some $\alpha_i\geq 0$. Now assume, as above, that ${\bf x}=(x_i)$ works for $\delta$
and ${\bf y}=(y_i)$ works for $\gamma$. By straight forward calculation we get that 
for each $i$,
\[
(r_j(\delta)+r_j(\gamma)-(x_i+y_i),r_j(\delta)+r_j(\gamma))>0
\]
if $i\neq j$, and 
\[
(r_i(\delta)+r_i(\gamma)-(x_i+y_i),r_i(\delta)+r_i(\gamma))=0.
\]
Taking into account that, for each $i$, $r_i(\delta+\gamma)=r_i(\delta)+r_i(\gamma)$,
the result follows. This can be checked using part 2 of 
Proposition~\ref{gettingcofdelta.prop}.
\end{proof}

\begin{definition}       \label{coterie.def}
We refer to $\mathscr{H}$ as the {\em coterie} of $\mathbb{H}$ or 
of $M_\mathbb{H}$.
\end{definition}

\begin{remark}  \label{type.rk}
Let $M$ be a reductive  monoid. The {\bf type map}
$\lambda : \Lambda\to 2^S$ of $M$ is a certain map from the
set of $G\times G$-orbits of $M$ to the set of subsets of the simple
roots. This type map determines $M$, as an abstract monoid, to within a kind
of central extension. This is equivalent to specifying the {\bf colored
face lattice} $\Lambda\subseteq 2^S\times 2^{\Lambda^1}$ of $M$
considered as a spherical variety (see Proposition 5.20 of \cite{bible}).
The type map of $M_\mathbb{H}$ is completely determined by the arrangement
$\mathbb{H}$. Each type map of a semisimple monoid associated with $\mathbb{H}$
is a kind of ``realization" of the type map of $M_\mathbb{H}$.
\end{remark}

\begin{remark} \label{coterie.rk}
Let $\delta\in\mathscr{H}$ and let $X_\delta=(M\backslash\{0\})/K^*$. 
Then any ample Cartier divisor $H$ on
$X_\delta$ will allow one to recover $X_\delta$ from the graded algebra of
global sections $\oplus_{n\geq0}\Gamma(X_\delta,nH)$. This same graded algebra for an
``ample" {\em Weil} divisor $H$ on $X_\delta$ is naturally the cone on some other projective
variety $X_\gamma$, a ``morphed" version of $X_\delta$. The resulting
variety $X_\gamma$ will agree with $X_\delta$ outside a closed subset of
codimension two. $\mathscr{H}$ can be thought of as the set of all possible
$W$-invariant, rational polytopes that can be constructed using $\mathbb{H}$ 
as the set of oriented hyperplanes associated with the facets.

The set of type maps partitions the cone  $\mathscr{H}\subseteq Cl(X)$ into the 
disjoint union of potentially smaller ample Cartier cones
$\mathscr{H}_\lambda$, each associated with some type map $\lambda$.
Thus one may think of $\mathscr{H}$ as the cone of ``ample Weil divisors" on 
$X_\delta=(M_\delta\backslash\{0\})/K^*$. 
\end{remark}

\section{$\mathscr{H}$ for $Env(G_0)$} \label{thediva.sec}

In this section we calculate explicitly the coterie $\mathscr{H}$ for 
$Env(G_0)$, where $G_0$ is a simple group. We obtain also some important
information about the rational polyhedral cone $\overline{\mathscr{H}}$. However 
we first explain how $Env(G_0)$ is related to other important constructions.

Let $L$ be a semisimple group of adjoint type, and suppose
that $\sigma :L\to L$ is an involution (so that
$\sigma\circ\sigma =id_L$) with $H=\{x\in L\mid\sigma (x)=x\}$.  The
{\bf plongement magnifique} of $L/H$ is the unique normal $L$-equivariant 
compactification $X$ of $L/H$ obtained by considering an irreducible 
representation $\rho :L\to Gl(V)$ of $L$
with $\dim (V^H)=1$ and with highest weight in general
position.  Then let $h\in V^H$ be nonzero and define
\[ X = \overline{\rho {(K)}[h]}\subseteq\mathbb{P}(V), \]
the Zariski closure of the orbit of $[h]$.  (See Section 2 of 
\cite{DP} for details.)
We are here concerned only with the case where $L=G_0\times G_0$,
$\sigma(g,h)=(h,g)$ and $G_0$ is a simple group. This amounts to the 
situation where we consider two-sided compactifications of $G_0$, or {\bf 
group embeddings}.

The relationship with semisimple monoids is as follows. 

If $G_0$ is a simple group and $M$ is a semisimple monoid 
with unit group $G_0\times K^*$ then we say $M$ is a {\bf canonical monoid} 
if there exists a finite morphism 
\[
\rho : M\to End(V)
\]
of algebraic monoids such that 
\begin{enumerate}
	\item[a)] $\rho$ is irreducible considered as a representation of $G_0$.
	\item[b)]the highest weight $\lambda$ of $\rho$ is of the form
	      $\lambda=\sum_\alpha c_\alpha\lambda_\alpha$ where 
	      $\{\lambda_\alpha\}$ is the set of fundamental dominant 
	      weights, and each $c_\alpha$ is nonzero.
\end{enumerate}
Canonical monoids are discussed in detail in \cite{PR}. There are many 
interesting characterizations of this class of monoids. 
It follows from the results of \cite{R2.5} or \cite{DP} that if $G_0$
is also of adjoint type then $(M\backslash\{0\})/K^*$ is isomorphic to the plongement 
magnifique. 

But it is important for us to consider all simple groups, and not just 
adjoint groups, since geometric properties like smoothness are not 
involved in our calculations. Also our mission here is to calculate the 
rational cone $\mathscr{H}$ whose lattice points correspond to 
semisimple monoids with $W$-arrangement $\mathbb{H}_\Delta$. Here we define 
$\mathbb{H}_\Delta$, by setting
\[
\Lambda^1\mathbb{H}_\Delta=\{(Span_\mathbb{Z}(\Delta\backslash\alpha),
\nu_\alpha)\;|\;\alpha\in\Delta\},
\]
where $\nu_\alpha : P\to\mathbb{Z}$ is defined by $\nu_\alpha(\beta)=-\delta_{\alpha\beta}\alpha$, 
and $\delta_{\alpha\beta}$ is the Kronecker delta. By Theorem~\ref{mvisisotomh.thm} 
we obtain that 

\begin{theorem}   \label{divavsenv.thm}
Let $M$ be a semisimple monoid and let $\mathbb{H}=\mathbb{H}(M)$.
The following are equivalent. 
\begin{enumerate}
	\item $M_\mathbb{H}=Env(G_0)$.
	\item $\mathbb{H}(M)=\mathbb{H}_\Delta$.
	\item There is a finite dominant morphism $M_V\to Env(G_0)$.
\end{enumerate}
\end{theorem}

In particular, if $M$ is a canonical monoid, then $M_\mathbb{H}\cong Env(G_0)$. Indeed, one 
can check directly that, for a canonical monoid $M$, $\mathbb{H}(M)=\mathbb{H}_\Delta$.
Thus $Env(G_0)$ is also the $\mathbb{Q}$-factorial hull of the plongement magnifique. 
As we shall see, there are many other semisimple monoids that have the same 
$\mathbb{Q}$-factorial hull as a canonical monoid. There are many a manifique 
wannabe.

We now proceed to the calculation.

Because of the special nature of canonical monoids and $Env(G_0)$ we are able 
to identify $\mathscr{H}$ as a subset of $X(T_0)\otimes\mathbb{Q}$ containing
the interior of the (rational) Weyl chamber $\mathscr{C}^0$. This is convenient 
for our calculations, but it also depicts $\mathscr{H}$ as a kind of virtual 
Borel-Weil-Bott theorem for rank-one, locally free sheaves on 
$(M\backslash CD2)/K^*$, where $CD2\subseteq M$ is the union of all
$G\times G$-orbits of codimension two or more. The set of canonical monoids
correspond to the points of $\mathscr{C}^0\subseteq\mathscr{H}$.
Indeed, we can get an inclusion
\[
\mathscr{H}\subseteq P_0\subseteq X(T_0)\otimes\mathbb{Q}
\]
by observing that for each $\delta\in\mathscr{H}$ the cone $\mathscr{C}(\delta)$
is of the form

\[
\mathscr{C}(\delta)=\mathscr{C}\cap(y - P_0)
\]
for some unique $y=y(\delta)\in P_0$, where $P_0$ is the rational cone
generated by the set of positive roots. This is so because  
the set of codimension-one colored faces of $\mathscr{C}(\delta)$ 
determines a collection of $n$ hyperplanes ($n=dim(T_0)$) which 
intersect at exactly one point somewhere in $X(T_0)\otimes\mathbb{Q}$,
but not necessarily in $\mathscr{C}(\delta)$.

Notice that
$\mathscr{C}\subseteq P_0$. Thus, we obtain inclusions 

\[
\mathscr{C}^0\subseteq\mathscr{H}\subseteq P_0\subseteq
X(T_0)\otimes\mathbb{Q}.
\]
The second of these inclusions is defined, as above, by
\[
\delta\leadsto y(\delta).
\]

Associated with the simple group $G_0$ is its {\bf Cartan matrix} $C$. The columns of 
$(C^{-1})^T$ are the coefficients that express the fundamental weights in terms of 
the simple roots. We refer to {\em Table 2} of \cite{OV} for the complete list of 
these coefficients.

Let $C_0(\Delta\backslash\{\alpha\})^0=
\{\sum_{\beta\neq\alpha}a_\beta\beta\in P_0\;|\; a_\beta > 0\;\text{for all}\;\beta\in\Delta\backslash\{\alpha\}\}$
be the interior of the rational cone
$C_0(\Delta\backslash\{\alpha\})$.

\begin{theorem}   \label{Hforthediva.thm}
Let $X$ be the plongement magnifique for the simple group $G_0$.
Let $P_0^0\subseteq P_0$ be the interior of the cone $P_0$ (i.e. 
$x=\sum a_\alpha\alpha\in P_0$ such that $a_\alpha > 0$ for all $\alpha$). 
Then in terms of the above
identification $\mathscr{H}\subseteq P$, 
\[
\mathscr{H}=\{x\in P^0\;|\;\text{for all $\alpha\in\Delta$ there exists 
$r_\alpha(x)>0$ such that}\; x-r_\alpha(x)\lambda_\alpha\in C_0(\Delta\backslash\{\alpha\})^0\}.
\]
If we write $x=\sum_{\beta\in\Delta}a_\beta\beta\in P^0_0$ then the following are equivalent.
\begin{enumerate}
\item $x\in\mathscr{H}$.
\item 
 \begin{enumerate}
 \item $a_\beta > \frac{c_{\alpha,\beta}}{c_{\alpha,\alpha}}a_\alpha$ for all $\alpha\neq\beta$,
 \item $a_\alpha > 0$ for all $\alpha\in\Delta$.
 \end{enumerate}
 \item
 \begin{enumerate}
 \item $a_\beta > \frac{c_{\alpha,\beta}}{c_{\alpha,\alpha}}a_\alpha$ 
 whenever $s_\alpha s_\beta\neq s_\beta s_\alpha$.
 \item $a_\alpha > 0$ for all $\alpha\in\Delta$.
 \end{enumerate}
\end{enumerate}
\end{theorem}
\begin{proof}
By Lemma 7.2.2 of \cite{bible} there is a unique $e_\alpha\in\Lambda$
such that $\lambda(e_\alpha)=\Delta\backslash\alpha$. Furthermore, by
Lemma 7.2.5 of \cite{bible} each $e_\alpha$ is in $\Lambda^1$ since
\[
F_\alpha = (z_\alpha + H_\alpha)\cap\mathscr{P}(\delta)
\]
is a codimension one face of $\mathscr{P}(\delta)$. Now each $e_\alpha$
has the property that
\[
\mathbb{Q}^+\lambda_\alpha\cap F_\alpha = \{r_\alpha\lambda_\alpha\}\in F_\alpha^0
\]
for some unique $r_\alpha\in\mathbb{Q}^+$, since
$W_{\Delta\backslash\alpha}$ acts on $F_\alpha$ with
$\{r_\alpha\lambda_\alpha\}$ as its unique fixed point.

But $\{r_\alpha\lambda_\alpha\}\in x-C_0(\Delta\backslash\{\alpha\})$ also,
and thus,
\[
\{r_\alpha\lambda_\alpha\}\in (x-C_0(\Delta\backslash\{\alpha\}))\cap
F_\alpha^0
\subseteq x-C_0(\Delta\backslash\{\alpha\})^0.
\]

Now let $x\in\mathscr{H}\subseteq P^0$. Then (as above) for each $\alpha\in\Delta$
there exists $r_\alpha(x)\in\mathbb{Q}$ such that

\[
x-r_\alpha(x)\lambda_\alpha\in C_0(\Delta\backslash\{\alpha\})^0.
\]
Write $x=\sum_{\beta\in\Delta}a_\beta\beta$ where
$(C^T)^{-1}=(c_{\beta,\alpha})$ and $C$ is the Cartan matrix. Thus

\[
x-r_\alpha(x)\lambda_\alpha=\sum_\beta(a_\beta-r_\alpha(x)c_{\beta,\alpha})\beta
\]
so that $a_\alpha-r_\alpha(x)c_{\alpha,\alpha}=0$. Thus,
\[
r_\alpha(x)=\frac{a_\alpha}{c_{\alpha,\alpha}}.
\]
If $\beta\neq\alpha$ we obtain that
\[
a_\beta-r_\alpha(x)c_{\beta,\alpha} > 0,
\]
so that
\[
a_\beta > \frac{c_{\beta,\alpha}}{c_{\alpha,\alpha}}a_\alpha
\]
for all $\beta\neq\alpha$.
We conclude that
\[
\mathscr{H}=\{x=\sum a_\beta\beta\;|\;a_\beta > 0\;\text{and}\;a_\beta >
\frac{c_{\beta,\alpha}}{c_{\alpha,\alpha}}a_\alpha\;\text{for all}\;\beta\neq\alpha\}
\]
It is sufficient to impose the condition ``$a_\beta >
\frac{c_{\beta,\alpha}}{c_{\alpha,\alpha}}a_\alpha$" in cases where
$s_\alpha s_\beta\neq s_\beta s_\alpha$. 
This follows from the fact that if we have 
\[
\alpha-\beta-....-\gamma
\]
representing a subdiagram of the Dynkin diagram, then
\[
c_{\alpha,\gamma}=\frac{c_{\alpha,\beta}}{c_{\beta,\beta}}c_{\beta,\gamma}.
\]
This is easily checked by inspecting each matrix $(C^T)^{-1}$. See Table 2
of \cite{OV}.
\end{proof}

\subsection{The Calculation of $\mathscr{H}$ for $X$ in each Case}

Using Theorem~\ref{Hforthediva.thm} we now calculate the cone $\mathscr{H}$ for the
plongement manifique associated with
each simple group $G_0$. In each case, we use $a_j$ instead of $a_{\alpha_j}$, 
with the numbering as dictated by the associated Dynkin diagram.

\vskip 14mm
\setlength{\unitlength}{1mm}
\begin{picture}(60,20)
\put(0, 23) {{\bf \underline{$A_n$}}} 

\put(15,23) {\circle{2.5}}

\put (24.5,23){\circle{2.5}}

\put (48.5,23) {\circle{2.5}}

\put (58,23) {\circle{2.5}}

\put (14,26) {{\scriptsize $1$}}

\put (24,26) {{\scriptsize $2$}}

\put (45,26) {{\scriptsize $n-1$}}

\put (57,26) {{\scriptsize $n$}}

\put (16, 23) {\line (1,  0) {6.8}}

\put (25.5, 23) {\line (1,  0) {4}}

\put (32, 22.7) {.\,.\,.\,.\,.\,.}

\put (47, 23) {\line (-1, 0) {4}}

\put (49.5, 23) {\line (1, 0) {6.8}}

\end{picture}

\vskip-18mm
$\mathscr{H}=\{x=\sum_ja_j\alpha_j\;|\;\text{$a_j$ as follows}\}$ 
\begin{enumerate}
\item $a_j > 0$.
\item $a_j > \frac{j}{j+1}a_{j+1}$ for $j<n$.
\item $a_j > \frac{n+1-j}{n+1-j+1}a_{j-1}$ for $j > 1$.
\end{enumerate}
----------------------------------------------------
\vskip 10mm

\begin{picture}(60,20)
\put (0, 20) {{\bf \underline{$B_n$}}} 

\put (15,20) {\circle{2.5}}

\put (24,20) {\circle{2.5}}

\put (49.5,20) {\circle{2.5}}

\put (58,20) {\circle{2.5}}

\put (14.1, 23) {{\scriptsize $1$}}

\put (23.6, 23) {{\scriptsize $2$}}

\put (45.5,23) {{\scriptsize $n-1$}}

\put (58,23) {{\scriptsize $n$}}

\put (16,20.1) {\line (1,0) {6}}

\put (31, 19.9) {.\,.\,.\,.\,.\,.\,.}

\put (25.2,20.1) {\line (1,0) {4}}

\put (44, 20.1) {\line (1,0) {4}}

\put (50.8, 20.5) {\line (1,0) {6}}

\put (50.8, 19.7) {\line (1,0) {6}}

\put (53, 19) {$>$}

\end{picture}

\vskip -14mm
$\mathscr{H}=\{x=\sum_ja_j\alpha_j\;|\;\text{$a_j$ as follows}\}$
\begin{enumerate}
\item $a_j > 0$.
\item $a_j > \frac{j}{j+1}a_{j+1}$ for $j<n$.
\item $a_j > a_{j-1}$ for $j > 1$.
\end{enumerate}
----------------------------------------------------
\vskip 10mm

\begin{picture}(60,15)

\put (0, 18) {{\bf \underline{$C_n$}}} 

\put (15,18) {\circle{2.5}}

\put (23.6,18) {\circle{2.5}}

\put (31, 17.6) {.\,.\,.\,.\,.\,.\,.}

\put (49.5,18) {\circle{2.5}}

\put (52.8,17) {$<$}

\put (58,18) {\circle{2.5}}

\put (14,21.2) {{\scriptsize $1$}}
\put (23.3,21.2) {{\scriptsize $2$}}

\put (45.8,21.2) {{\scriptsize $n-1$}}

\put (57.3,21.2) {{\scriptsize $n$}}

\put (16,18.2) {\line (1,0) {6}}

\put (25, 18.2) {\line (1, 0) {4}}

\put (44,18.2) {\line (1,0) {4}}

\put (50.8, 17.8) {\line (1,0){6}}

\put (50.8, 18.6) {\line (1,0){6}}

\end{picture}

\vskip-14mm
$\mathscr{H}=\{x=\sum_ja_j\alpha_j\;|\;\text{$a_j$ as follows}\}$
\begin{enumerate}
\item $a_j > 0$.
\item $a_j > a_{j-1}$ for $j<n$.
\item $a_j > \frac{j}{j+1}a_{j+1}$ for $j<n-1$.
\item $a_n > \frac{1}{2}a_{n-1}$.
\item $a_{n-1} > \frac{2(n-1)}{n}a_n$.
\end{enumerate}
----------------------------------------------------
\vskip 16mm

\begin{picture}(60,20)
\put (0, 25) {{\bf \underline{$D_n$}}} 

\put (14.8,25) {\circle{2.5}}

\put (23.6,25) {\circle{2.5}}

\put (49.5,25) {\circle{2.5}}

\put (58,30) {\circle{2.5}}

\put (58,20) {\circle{2.5}}

\put (14,28.5) {{\scriptsize $1$}}

\put (23,28.5) {{\scriptsize $2$}}

\put (45.5,28.5){{\scriptsize $n-2$}}

\put (54.5,33) {{\scriptsize $n-1$}}

\put (57,14.5) {{\scriptsize $n$}}

\put (16, 25) {\line (1,  0) {6}}      

\put (25, 25) {\line (1,  0) {4}}

\put (31, 24.5) {.\,.\,.\,.\,.\,.\,.}

\put (44, 25) {\line (1,0) {4}}

\put (50.6, 25.2) {\line (2,1) {7}}

\put (50.6, 24.4) {\line (2,  -1) {7}}

\end{picture}

\vskip-14mm
$\mathscr{H}=\{x=\sum_ja_j\alpha_j\;|\;\text{$a_j$ as follows}\}$
\begin{enumerate}
\item $a_j > 0$.
\item $a_j > \frac{j}{j+1}a_{j+1}$ for $j<n-2$.
\item $a_{n-2} > \frac{2(n-2)}{n}a_{n-1}$
\item $a_{n-2} > \frac{2(n-2)}{n}a_{n}$
\item $a_n > \frac{1}{2}a_{n-2}$
\item $a_{n-1} > \frac{1}{2}a_{n-2}$
\end{enumerate}
----------------------------------------------------
\vskip 18mm

\setlength{\unitlength}{1mm}
\begin{picture}(80,20)
\put (0, 30) {{\bf \underline{$E_6$}}}  

\put (14.8,30) {\circle{2.5}}

\put (23.5,30) {\circle{2.5}}

\put (32.1,30) {\circle{2.5}}

\put (41,30){\circle{2.5}}

\put (49.6,30) {\circle{2.5}}

\put (32.1,23.6) {\circle{2.5}}

\put (13.9,32.8) {{\scriptsize $1$}}

\put (22.5,32.8) {{\scriptsize $2$}}

\put(31.2,32.8) {{\scriptsize $3$}}

\put (40.3,32.8) {{\scriptsize $4$}}

\put (49,32.8) {{\scriptsize $5$}}

\put (35,22.3){{\scriptsize $6$}}

\put (16, 30.1) {\line (1,0) {6}}

\put (24.7, 30.1) {\line (1,0) {6}}

\put (33.5, 30.1) {\line (1,0) {6}}

\put (42.2, 30.1) {\line (1,0) {6}}

\put (32.1, 29) {\line (0, -1) {4}}

\end{picture}

\vskip-14mm
$\mathscr{H}=\{x=\sum_ja_j\alpha_j\;|\;\text{$a_j$ as follows}\}$
\begin{enumerate}
\item $a_j > 0$.
\item $8a_1 > 4a_2 > 5a_1$.
\item $5a_3 > 6a_2 > 4a_3$.
\item $5a_3 > 6a_4 > 4a_3$.
\item $4a_3 > 6a_6 > 3a_3$.
\item $8a_5 > 4a_4 > 5a_5$.
\end{enumerate}
----------------------------------------------------
\vskip 14mm

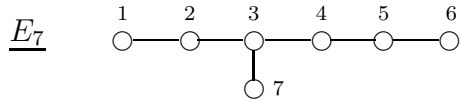
\begin{picture}(60,20)
\put (0, 23) {{\bf \underline{$E_7$}}} 

\put (15,23) {\circle{2.5}} 

\put (24,23) {\circle{2.5}}

\put(32.5,23){\circle{2.5}}

\put(41.5,23) {\circle{2.5}}

\put (49.7,23){\circle{2.5}}

\put (58.5,23) {\circle{2.5}}

\put (32.5,16.5) {\circle{2.5}}

\put (14.3,25.7) {{\scriptsize $1$}}

\put (23.2,25.7) {{\scriptsize $2$}}

\put (31.7,25.7) {{\scriptsize $3$}}

\put(40.7,25.7) {{\scriptsize $4$}}

\put (49,25.7) {{\scriptsize $5$}}

\put (58,25.7) {{\scriptsize $6$}}

\put (35,15.8){{\scriptsize $7$}}

\put(16.4, 23) {\line (1,0) {6}}

\put (25, 23) {\line (1,0) {6}}

\put (34, 23) {\line (1,0) {6}}

\put (42.6, 23) {\line (1,  0){6}}

\put (51, 23) {\line (1,  0) {6}}

\put (32.4, 21.7) {\line (0,-1) {4}}

\end{picture}

\vskip-14mm
$\mathscr{H}=\{x=\sum_ja_j\alpha_j\;|\;\text{$a_j$ as follows}\}$
\begin{enumerate}
\item $a_j > 0$.
\item $3a_1 > 4a_2 > 2a_1$.
\item $6a_2 > 4a_3 > 5a_2$.
\item $10a_4 > 12a_3 > 9a_4$.
\item $9a_5 > 6a_4 > 8a_5$.
\item $7a_4 > 12a_7 > 6a_4$.
\item $4a_6 > 2a_5 > 3a_6$.
\end{enumerate}
----------------------------------------------------
\vskip 14mm

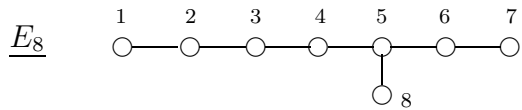
\begin{picture}(60,16)
\put (0,23) {{\bf \underline{$E_8$}}} 

\put(15,23){\circle{2.5}}

\put (24,23) {\circle{2.5}}

\put (32.7,23) {\circle{2.5}}

\put (41,23){\circle{2.5}}

\put (49.5,23) {\circle{2.5}}

\put (58,23) {\circle{2.5}}

\put(66.5,23) {\circle{2.5}}

\put (49.5,16.6) {\circle{2.5}}

\put (14.1,26.2){{\scriptsize $1$}}

\put (23.3,26.2) {{\scriptsize $2$}}

\put (32,26.2) {{\scriptsize $3$}}

\put (40.6,26.2) {{\scriptsize $4$}}

\put (48.7,26.2) {{\scriptsize $5$}}

\put (57.1,26.2) {{\scriptsize $6$}}

\put (66,26.2) {{\scriptsize $7$}}

\put (52,15) {{\scriptsize $8$}}

\put (16.3, 23) {\line (1,0) {6}}

\put (25.3, 23) {\line (1,0) {6}}

\put (34,23) {\line (1,0) {6}}

\put (42.4,23) {\line (1,0) {6}}

\put (50.7,23) {\line (1,0) {6}}

\put (59.3,23) {\line (1,0) {6}}

\put (49.5, 22) {\line (0,-1) {4}}

\end{picture}

\vskip-14mm
$\mathscr{H}=\{x=\sum_ja_j\alpha_j\;|\;\text{$a_j$ as follows}\}$
\begin{enumerate}
\item $a_j > 0$.
\item $4a_1 > 2a_2 > 3a_1$.
\item $9a_2 > 6a_3 > 8a_2$.
\item $16a_3 > 12a_4 > 15a_3$.
\item $25a_4 > 20a_5 > 24a_4$.
\item $16a_8 > 8a_5 > 15a_8$.
\item $21a_6 > 14a_5 > 20a_6$.
\item $8a_7 > 4a_6 > 7a_7$.
\end{enumerate}
----------------------------------------------------
\vskip 14mm

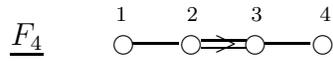
\begin{picture}(60,11)
\put (0, 14) {{\bf \underline{$F_4$}}} 

\put (15,14) {\circle{2.5}}

\put(24,14) {\circle{2.5}}

\put (27,13.1) {$>$}

\put (32.5,14) {\circle{2.5}}

\put (41.5,14) {\circle{2.5}}

\put (14,17.2) {{\scriptsize $1$}}

\put (23.4,17.2) {{\scriptsize $2$}}

\put (32,17.2) {{\scriptsize $3$}}

\put(41.1,17.2) {{\scriptsize $4$}}

\put (16.3,14.2) {\line (1,0){6}}

\put (25.3,14.6) {\line (1,0){6}}

\put (25.3,13.6) {\line (1,0){6}}

\put (33.8,14.2) {\line (1,0){6}}

\end{picture}

\vskip-10mm
$\mathscr{H}=\{x=\sum_ja_j\alpha_j\;|\;\text{$a_j$ as follows}\}$
\begin{enumerate}
\item $a_j > 0$.
\item $4a_1 > 2a_2 > 3a_1$.
\item $9a_2 > 12a_3 > 8a_2$.
\item $4a_3 > 6a_4 > 3a_3$.
\end{enumerate}
----------------------------------------------------
\vskip 14mm

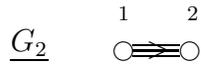
\begin{picture}(60,25)

\put (0,27) {{\bf \underline{$G_2$}}}  

\put (15,27){\circle{2.5}}

\put (18,26) {$>$}

\put (23.8,27){\circle{2.5}}

\put(14.3,31) {{\scriptsize $1$}}

\put (23.5,31) {{\scriptsize $2$}}

\put (16.3,27.7) {\line (1,0){6.1}}

\put (16.3,27) {\line (1,0){6}}

\put (16.3,26.3) {\line (1,0){6.1}}

\end{picture}

\vskip-14mm
$\mathscr{H}=\{x=\sum_ja_j\alpha_j\;|\;\text{$a_j$ as follows}\}$
\begin{enumerate}
\item $a_j > 0$.
\item $4a_2 > 2a_1 > 3a_2$. 
\end{enumerate}
\vskip 14mm

Notice that, for any rank-two example, $\mathscr{H}=\mathscr{C}$ the Weyl chamber.

\subsection{The Structure of ${\mathscr{H}}$ and $\overline{\mathscr{H}}$}

$\overline{\mathscr{H}}$ is a rational polyhedral cone. It is of interest to 
identify the face lattice $\mathscr{F}$ of $\overline{\mathscr{H}}$.
It turns out that the faces of 
$\overline{\mathscr{H}}$ are indexed by orientations of the associated
Dynkin diagram $\mathscr{D}$. To describe this correspondence let 
\[
\mathscr{E}=\{(\alpha,\beta)\in\Delta\times\Delta\;|\;s_\alpha s_\beta\neq s_\beta s_\alpha\}.
\]
be the set of edges of $\mathscr{D}$ and let 
\[
\Gamma = \{\leftarrow,--,\rightarrow\}=\{l,n,r\}
\]
be the realm of possible orientations (left, neutral, right) of each edge. Define
\[
F=\{f\in Hom(\mathscr{E},\Gamma)\;|\; f(\alpha,\beta)=r\;\text{if and only if}\;f(\beta,\alpha)=l\}.
\]
If $f\in F$ then, for each $(\alpha,\beta)\in\mathscr{E}$, $f(\alpha,\beta)$ 
represents either an arrow from $\alpha$ to $\beta$ ($f(\alpha,\beta)=r$ and
$f(\beta,\alpha)=l$), a broken line between $\alpha$ and $\beta$ 
($f(\alpha,\beta)=f(\beta,\alpha)=n$), or an arrow from $\beta$ to $\alpha$
($f(\alpha,\beta)=l$ and $f(\beta,\alpha)=r$).
So we can think of these faces as diagrams like the following.
\[
\alpha_1\leftarrow\alpha_2 --\alpha_3\rightarrow\alpha_4\leftarrow\alpha_5
\]
We now define the ordering on $F$. If $f,g\in F$ we define 
$f\geq g$ if 
\begin{enumerate}
\item[a)] $g(\alpha,\beta)=n$ implies that $f(\alpha,\beta)=n$,
\item[b)] $f(\alpha,\beta)=r$ implies that $g(\alpha,\beta)=r$ and
\item[c)] $f(\alpha,\beta)=l$ implies that $g(\alpha,\beta)=l$.
\end{enumerate}
Thus $f\geq g$ if $f$ can be obtained from $g$ by replacing some of its 
arrows by broken lines. It is easy to check that $(F,\geq)$ is isomorphic 
to the face lattice of a cube. The vertices of $F$ are those orientations $f\in F$
such that, for all $(\alpha,\beta)\in\mathscr{E}$, $f(\alpha,\beta)\neq n$.

\begin{theorem}   \label{fisacube.thm}
Let $G_0$ be a simple group of rank $n$. There is a canonical one-to-one 
correspondence
\[
(F,\geq)\cong\mathscr{F}.
\]
In particular, the vertex figure of $\overline{\mathscr{H}}$ is isomorphic to 
the face lattice of an $(n-1)$-cube.
\end{theorem}
\begin{proof}
In each case of rank $n$, $\mathscr{H}$ is defined by $n-1$ conditions of the form
\[
ra_i > sa_j > ta_i.
\]
Furthermore, the edges of $\mathscr{D}$ are exactly the pairs $(i,j)$ that occur 
in each list of conditions. Now $\overline{\mathscr{H}}$ is defined by the $n-1$ 
conditions
\[
ra_i \geq sa_j \geq ta_i.\;\;\;\;\;\;\;(*)
\]
A face of $\overline{\mathscr{H}}$ is defined by replacing at most one ``$\geq$" in 
each of these conditions by an ``=". Replacing the left ``$\geq$" by a ``=" results in  
a $i\rightarrow j$. Replacing the right ``$\geq$" by a ``=" results in  a $i\leftarrow j$. 
Finally, if no ``$\geq$" in $(*)$ is replaced, then the result is a $i -- j$. This determines 
the canonical isomorphism $(F,\geq)\cong\mathscr{F}$ of partially ordered sets.
\end{proof}

\begin{example} \label{extremalfortypea4.ex}
In this example we calculate the extremal rays of $\overline{\mathscr{H}}$ 
in case $G_0$ is a simple group of type $A_4$. In each case we include the 
associated oriented Dynkin diagram, a representative of the corresponding 
extremal ray and the defining conditions for this extremal ray. It follows 
from Theorem~\ref{fisacube.thm} above that each of these extremal rays 
corresponds to an oriented diagram $f\in F$ such that $f(\alpha,\beta)\neq n$ for all 
$(\alpha,\beta)\in\mathscr{E}$.\\

\begin{enumerate}
\item $\alpha_1\rightarrow\alpha_2\rightarrow\alpha_3\rightarrow\alpha_4$
 \begin{enumerate}
  \item $(1/4,1/2,3/4,1)$
  \item $2a_1=a_2, 3a_2=2a_3, 4a_3=3a_4.$
 \end{enumerate}
\item $\alpha_1\rightarrow\alpha_2\rightarrow\alpha_3\leftarrow\alpha_4$
 \begin{enumerate}
  \item $(2/3,4/3,2,1)$
  \item  $2a_1=a_2, 3a_2=2a_3, a_3=2a_4.$
 \end{enumerate}
\item $\alpha_1\rightarrow\alpha_2\leftarrow\alpha_3\rightarrow\alpha_4$
 \begin{enumerate}
  \item $(9/16,9/8,3/4,1)$
  \item $2a_1=a_2, 2a_2=3a_3, 4a_3=3a_4.$
 \end{enumerate}
\item $\alpha_1\rightarrow\alpha_2\leftarrow\alpha_3\leftarrow\alpha_4$
 \begin{enumerate}
  \item $(3/2,3,2,1)$
  \item $2a_1=a_2, 2a_2=3a_3, a_3=2a_4.$
 \end{enumerate}
\item $\alpha_1\leftarrow\alpha_2\rightarrow\alpha_3\rightarrow\alpha_4$
 \begin{enumerate}
  \item $(2/3,1/2,3/4,1)$
  \item $3a_1=4a_2, 3a_2=2a_3, a_3=2a_4.$
 \end{enumerate}
\item $\alpha_1\leftarrow\alpha_2\rightarrow\alpha_3\leftarrow\alpha_4$ 
 \begin{enumerate}
  \item $(16/9,4/3,2,1)$
  \item $3a_1=4a_2, 3a_2=2a_3, a_3=2a_4.$
 \end{enumerate}
\item $\alpha_1\leftarrow\alpha_2\leftarrow\alpha_3\rightarrow\alpha_4$
 \begin{enumerate}
  \item $(3/2,9/8,3/4,1)$
  \item $3a_1=4a_2, 2a_2=3a_3, 4a_3=3a_4.$
 \end{enumerate}
\item $\alpha_1\leftarrow\alpha_2\leftarrow\alpha_3\leftarrow\alpha_4$
 \begin{enumerate}
  \item $(4,3,2,1)$
  \item $3a_1=4a_2, 2a_2=3a_3, a_3=2a_4.$
 \end{enumerate}
\end{enumerate}
\end{example}

\vspace{20pt}

\noindent Lex E. Renner \\
Department of Mathematics \\
University of Western Ontario \\
London, N6A 5B7, Canada \\

\enddocument